\documentclass[fleqn]{mat01}
\usepackage{times,mathtimy,amssymb,latexsym}

\def\R{{\Bbb R}}

\newcommand{\Ques}[2]{}

\def\disp{\displaystyle}

\def\bc{\begin{center}}       \def\ec{\end{center}}
\def\ba{\begin{array}}        \def\ea{\end{array}}
\def\be{\begin{equation}}     \def\ee{\end{equation}}
\def\bea{\begin{eqnarray}}    \def\eea{\end{eqnarray}}
\def\beaa{\begin{eqnarray*}}  \def\eeaa{\end{eqnarray*}}

\def\lb{\label}               \def\x#1{(\ref{#1})}
\def\bt{\bibitem}

\def\oo{\infty}

\def\nn{\nonumber}

\begin{document}

\setcounter{page}{85} \firstpage{85}

\renewcommand\theequation{\thesection\arabic{equation}}

\newtheorem{theore}{Theorem}
\renewcommand\thetheore{\arabic{section}.\arabic{theore}}
\newtheorem{theor}[theore]{\bf Theorem}
\newtheorem{lem}[theore]{\it Lemma}
\newtheorem{propo}[theore]{\rm PROPOSITION}
\newtheorem{coro}[theore]{\rm COROLLARY}
\newtheorem{definit}[theore]{\rm DEFINITION}
\newtheorem{exa}[theore]{\it Example}
\newtheorem{rem}[theore]{\it Remark}

\newtheorem{theoree}{Theorem}
\renewcommand\thetheoree{\arabic{theoree}}
\newtheorem{case}[theoree]{\it Case}

\title{Positive solutions and eigenvalue intervals for nonlinear
systems}

\markboth{Jifeng Chu, Donal O'Regan and Meirong Zhang}{Existence
of positive solutions for nonlinear systems}

\author{JIFENG CHU$^{1,2}$, DONAL O'REGAN$^3$ and
MEIRONG ZHANG$^{1}$}

\address{$^{1}$Department of Mathematical Sciences,
Tsinghua University, Beijing~100~084, China\\
\noindent $^2$Department of Applied Mathematics, Hohai University,
Nanjing~210~098, China\\
\noindent $^3$Department of Mathematics, National University of
Ireland, Galway, Ireland\\
\noindent E-mail: jifengchu@yahoo.com.cn;
mzhang@math.tsinghua.edu.cn}

\volume{117}

\mon{February}

\parts{1}

\pubyear{2007}

\Date{MS received 22 February 2006}

\begin{abstract}
This paper deals with the existence of positive solutions for the
nonlinear system
\begin{align*}
(q(t)\phi(p(t)u'_{i}(t)))'+f^{i}(t,\textbf{u})=0,\quad 0<t<1,\quad
i=1,2,\dots,n.
\end{align*}
This system often arises in the study of positive radial solutions
of nonlinear elliptic system. Here
$\textbf{u}=(u_{1},\dots,u_{n})$ and $f^{i},~i=1,2,\dots,n$ are
continuous and nonnegative functions, $p(t),~q(t)\hbox{\rm :}\
[0,1]\rightarrow (0,\oo)$ are continuous functions. Moreover, we
characterize the eigenvalue intervals for
\begin{align*}
(q(t)\phi(p(t)u'_{i}(t)))'+\lambda h_{i}(t)g^{i} (\textbf{u})=0,
\quad 0<t<1,\quad i=1,2,\dots,n.
\end{align*}
The proof is based on a well-known fixed point theorem in cones.
\end{abstract}


\maketitle

\section{Introduction}

\setcounter{equation}{0}

In this paper we study the existence of positive solutions for the
nonlinear system
\begin{equation}\lb{a}
\left\{\begin{array}{ll}
(q(t)\phi(p(t)u'_{1}(t)))'+f^{1}(t,\textbf{u})=0,\quad 0<t<1,\\[.1pc]
\dots\\[.2pc]
(q(t)\phi(p(t) u'_{n}(t)))'+f^{n}(t,\textbf{u})=0,\quad
0<t<1,
\end{array}\right.
\end{equation}
with the following boundary condition
\begin{equation}
\lb{b}\textbf{u}(0)=0,\quad \textbf{u}(1)=0.
\end{equation}
Here $\phi(x)=|x|^{p-2}x,~p>1$, $\textbf{u}=(u_{1},\dots,u_{n})$.
We always make the following assumptions:

\begin{itemize} \leftskip .2pc
\item[$({\textmd{H}_{1}})$]
$f^{i}\hbox{\rm :}\ [0,1]\times\R^{n}_{+} \rightarrow(0,\oo)$ is
continuous, $i=1,2,\dots,n.$

\item[$({\textmd{H}_{2}})$] $p(t),~q(t)\hbox{\rm :}\ [0,1]\rightarrow (0,\oo)$
are continuous functions and $q(t)$ is nondecreasing on [0,1].
\end{itemize}

Problems \x{a} and \x{b} often arise from the study of positive
radial solutions for the nonlinear elliptic system of the form
\begin{equation}\lb{c}
\left\{\begin{array}{ll}
\textmd{div}(|\nabla u_{1}|^{p-2}\nabla
u_{1})+k_{1}(|x|)g^{1}(\textbf{u})=0,\\[.1pc]
\dots\\[.2pc]
\textmd{div}(|\nabla u_{n}|^{p-2}\nabla
u_{n})+k_{n}(|x|)g^{n}(\textbf{u})=0,\end{array}\right.
\end{equation}
in the domain $0<R_{1}<|x|<R_{2}<\oo,~x\in\R^{N},~N\geq 2$ with
the following boundary condition:
\begin{equation}
\lb{d}u_{i}=0\quad \mbox{on~}|x|=R_{1} \quad \mbox{and} \quad
|x|=R_{2},~i=1,2,\dots,n.
\end{equation}

In recent years, positive radial solutions for nonlinear elliptic
equation or elliptic system have been studied by many authors and
we refer the reader to \cite{b1,d1,d2,l1,l2,l3,l4}. It was proved
in \cite{b1,l4} that the classical elliptic equation
\begin{equation*}
\Delta u+\lambda k(|x|)f(u)=0,~\mbox{in}~R_{1}<|x| <R_{2},~x\in
\R^{N},~ N\geq 2
\end{equation*}
has at least one positive radial solution under the assumption
that $f$ is either superlinear or sublinear. The one-dimensional
$p$-Laplacian boundary value problem has also attracted
considerable attention \cite{a2,b2,c2,k2,o,w2}. In \cite{h1}, it
was proved that
\begin{equation*}
(\phi(u'))'+\lambda h(t)f(u)=0
\end{equation*}
with boundary condition $u(0)=u(1)=0$ having at least one positive
solution for certain finite intervals of $\lambda$ if one of
$f_{0}$ and $f_{\oo}$ is large enough and the other one is small
enough.

In this paper, we also study the following eigenvalue problem
\begin{equation}\lb{i}
\left\{\begin{array}{ll}
(q(t)\phi(p(t)u'_{1}(t)))'+\lambda
h_{1}(t)g^{1}(\textbf{u})=0,\quad 0<t<1,\\[.1pc]
\dots\\[.2pc]
(q(t)\phi(p(t)u'_{n}(t)))'+\lambda
h_{n}(t)g^{n}(\textbf{u})=0,\quad
0<t<1.\end{array}\right.
\end{equation}
We prove that \x{i} and \x{b} have at least one positive solution
for each $\lambda$ in an explicit eigenvalue interval. Recently,
several eigenvalue characterizations for different kinds of
boundary value problems have appeared and we refer the reader to
\cite{a1,a2,c1,c3,h1,h2}. In this paper, we will show how our
method allows to improve the range of eigenvalue intervals. The
new results are easily derived from a general result, stated as
Theorem~\ref{main}, which gives sufficient conditions to guarantee
the existence of at least one positive solution for systems \x{a}
and \x{b}.

Our arguments are based on a well-known fixed point theorem in
cones. Many authors \cite{a2,c1,f} have used this fixed point
theorem to discuss the existence of positive solutions for
different boundary value problems. In \cite{h2,w1}, the existence,
multiplicity and nonexistence of positive solutions for nonlinear
systems of ordinary differential equations were considered using
fixed point index. The same problems for the quasilinear elliptic
system \x{c} were studied in \cite{o2}. The novelty of this paper
is a choice of a cone different from that used in \cite{h2,o2,w1}.
Such a choice of cones (see \x{j}) is well fit with the systems we
are considering.

Finally, it is worth remarking here that, we can also deal with
the nonlinear system \x{a} with one of the following two sets of
boundary conditions
\begin{align*}
\rm{\textbf{u}}'(0)&=0,\qquad \rm{\textbf{u}}(1)=0,\\[.2pc]
\rm{\textbf{u}}(0)&=0,\qquad \rm{\textbf{u}}'(1)=0.
\end{align*}
However since the arguments are essentially the same (in fact
easier), we will restrict our discussion to boundary data \x{b}.

The notation used is as follows: $\R_{+}=[0,\oo),
\R^{n}_{+}=\prod_{i=1}^{n}\R_{+}.$ For $\textbf{u}=
(u_{1},\dots,u_{n}) \in\R^{n}_{+},~\|\textbf{u}\|=\max_{i=1,
2,\dots, n}|u_{i}|.$

The remaining part of the paper is organized as follows. In \S2,
some preliminary results are given and in \S3, the main results
are proved.

\section{Preliminaries}

\setcounter{equation}{0}

The proof of the main results is based on a well-known fixed point
theorem in cones. We recall the statement of this result below,
after introducing the definition of a cone.

\begin{definit}$\left.\right.$\vspace{.5pc}

\noindent {\rm Let $X$ be a Banach space and $K$ be a closed,
nonempty subset of $X$. $K$ is a cone if
\begin{enumerate}
\renewcommand\labelenumi{(\roman{enumi})}
\leftskip .15pc
\item $\alpha u+\beta v\in K$ for all
$u,v\in K$ and all $\alpha,\beta>0$,

\item $u, -u \in K$ implies $u=0.$\vspace{-.5pc}
\end{enumerate}}
\end{definit}

We also recall that a completely continuous operator means a
continuous operator which transforms every bounded set into a
relatively compact set. If $D$ is a subset $X$, we write
$D_{K}=D\cap K$ and $\partial_{K}D=(\partial D)\cap K.$

\begin{theor}[{\cite{k1}}]\lb{cone}
Let $X$ be a Banach space and $K\ (\subset X)$ be a cone. Assume
that $\Omega^1,\ \Omega^2$ are open subsets of $X$ with
$\Omega^1_{K}\neq\emptyset, \overline{\Omega^1}_{K} \subset
\Omega^2_{K}.$ Let
\begin{equation*}
T\hbox{\rm :}\ \overline{\Omega^2}_{K} \rightarrow K
\end{equation*}
be a continuous{\rm ,} completely continuous operator such that
either
\begin{enumerate}
\renewcommand\labelenumi{\rm (\roman{enumi})}
\leftskip .15pc
\item $\| Tu \| \geq \| u \|, u\in
\partial_{K}\Omega^{1}$ and $ \| Tu \| \leq \| u \|,u\in
\partial_{K}\Omega^{2}${\rm ;} or

\item $\| Tu \| \leq \| u \|, u\in \partial_{K}\Omega^{1}$ and $\|
Tu \| \geq \| u \|, u\in \partial_{K}\Omega^{2}$.\vspace{-.5pc}
\end{enumerate}
Then $T$ has a fixed point in $\overline{\Omega^2}_{K}\!\backslash
\Omega^1_{K}.$
\end{theor}

As usual, we denote by $C[0,1]$ the space of continuous functions
from [0,1] to $\R$. In $C[0,1]$ we shall consider the norm
$|u|_{0}=\max_{0\leq t\leq 1}|u(t)|$. In order to apply
Theorem~\ref{cone} below, we take $X=C[0,1]\times
C[0,1]\times\dots\times C[0,1]$ ($n$ times) with the norm
$\|\textbf{u}\|=\max_{i=1, 2, \dots,n}|u_{i}|_{0}$ for
$\textbf{u}=(u_{1},\dots, u_{n})\in X$. Then $X$ is a Banach
space. Define
\begin{equation}\lb{j}
\hskip -4pc K= \left\{\textbf{u}\in X\hbox{\rm :}\ u_{i}(t)\geq
0~, \forall ~t\in[0,1]~\mbox{and}~\min_{\frac{1}{4}\leq t\leq
\frac{3}{4}} u_{i}(t)\geq \rho|u_{i}|_{0},~i=1, 2, \dots,
n\right\},
\end{equation}
where $\rho$ is given by
\begin{equation}\lb{p}
\rho= \left[ \int_{0}^{1} \frac{1}{p(s)}{\rm d}s \right]^{-1} \min
\left\{ \int_{0}^{\frac{1}{4}}\frac{1} {p(s)}{\rm
d}s,~\int^{1}_{\frac{3}{4}}\frac{1}{p(s)}{\rm d}s\right\}.
\end{equation}
One can easily verify that $K$ is a cone in $X$.\pagebreak

Let $T\hbox{\rm :}\ K\rightarrow X$ be a map with components
$(T^{1},\dots, T^{n}),$ where $T^{i},i=1,2,\dots,n$ is defined by
\begin{equation}\lb{k}
(T^{i}\textbf{u})(t)= \left\{\begin{array}{ll}
\int_{0}^{t}\frac{1}{p(s)}\phi^{-1}
\left(\frac{1}{q(s)}\int_{s}^{\sigma_{i}}f^{i}(\tau,\textbf{u}(\tau))
{\rm d}\tau\right){\rm d}s,\quad 0\leq t\leq \sigma_{i},\\[.1pc]
\dots\\[.2pc]
\int_{t}^{1}\frac{1}{p(s)}\phi^{-1}
\left(\frac{1}{q(s)}\int_{\sigma_{i}}^{s}f^{i}(\tau,\textbf{u}(\tau))
{\rm d}\tau\right){\rm d}s,\quad \sigma_{i}\leq t\leq
1;\end{array}\right.
\end{equation}
here $\sigma_{i}\in(0,1)$ is a solution of the equation
\begin{equation}
\lb{l}\Theta^{i}\textbf{u}(t)=0,\quad 0\leq t\leq 1
\end{equation}
and the map $\Theta^{i}\hbox{\rm :}\ K\rightarrow C[0,1]$ is
defined by
\begin{align*}
\Theta^{i}\textbf{u}(t) &=\int_{0}^{t}\frac{1}{p(s)}\phi^{-1}
\left(\frac{1}{q(s)}\int_{s}^{t}f^{i}(\tau,\textbf{u}(\tau))
{\rm d}\tau\right){\rm d}s\\[.5pc]
&\quad\, -\int_{t}^{1}\frac{1}{p(s)}\phi^{-1}
\left(\frac{1}{q(s)}\int_{t}^{s}f^{i}(\tau,\textbf{u}(\tau)) {\rm
d}\tau\right){\rm d}s,\quad 0<t<1.
\end{align*}

\begin{lem}\hskip -.3pc {\rm \cite{h2}.} \ \ Assume $({\rm
H}_1)$ and $({\rm H}_2)$ hold. Then, for any ${\rm{\textbf{u}}}\in
K$ and $i=1,2,\dots,n$, $\Theta^{i}{\rm{\textbf{u}}}(t)=0$ has at
least one solution in $(0,1)$ and $T$ is well-defined.
\end{lem}

\begin{lem}\lb{har}\hskip -.3pc {\rm \cite{h2}.} \ \ Assume $({\rm
H}_2)$ holds. Let $u\in C^{1}[0,1]$ be a nonnegative function and
$q(t)\phi(p(t)u'(t))$ is nonincreasing on $[0,1]$. Then
\begin{align*}
u(t) \geq \left[\int_{0}^{1} \frac{1}{p(s)}{\rm
d}s\right]^{-1}\min \left\{\int_{0}^{t}\frac{1} {p(s)}{\rm
d}s,~\int_{t}^{1}\frac{1}{p(s)}{\rm d}s\right\}|u|_{0},~0\leq
t\leq 1.
\end{align*}
In particular{\rm ,}
\begin{equation*}
\min_{\frac{1}{4}\leq t\leq \frac{3}{4}}u(t)\geq
\rho|u|_{0}.
\end{equation*}
\end{lem}

\begin{rem}{\rm If $p(t)=C$ for all $0\leq t\leq 1,$ then
$\rho=\frac{1}{4}$, and this was used in \rm{\cite{o2}}.}
\end{rem}

\begin{lem}\hskip -.3pc{\rm \cite{h2}.} \ \ Assume $({\rm H}_1)$
and $({\rm H}_2)$ hold. Then $T(K)\subset K$. Moreover{\rm ,} $T$
is continuous and completely continuous.
\end{lem}

\begin{proof}
Lemma~\ref{har} implies that $T(K)\subset K$. It is easy to see
that $T$ is continuous and completely continuous since $({\rm
H}_1)$ and $({\rm H}_2)$ hold. See \cite{w1} for a proof.\hfill
$\Box$
\end{proof}

\section{Main results}

\setcounter{equation}{0}

\setcounter{theore}{0}

In this section we establish the existence of positive solutions
for \x{a} and \x{b}. Moreover, we characterize the eigenvalues for
problems \x{i} and \x{b}.

For the given function $a\in C[0,1],$ let
\begin{align}\lb{m}
\gamma_{a}(t) &= \frac{\rho}{2}\left[\int_{\frac{1}{4}}^{t}
\frac{1}{p(s)}\phi^{-1} \left(\frac{1}{q(s)}\int_{s}^{t}a(\tau)
{\rm d}\tau\right){\rm d}s \right.\nonumber\\[.5pc]
&\quad\, \left. + \int^{\frac{3}{4}}_{t}\frac{1}{p(s)}
\phi^{-1}\left( \frac{1}{q(s)}\int^{s}_{t}a(\tau){\rm
d}\tau\right){\rm d}s\right],\quad \frac{1}{4}\leq t\leq
\frac{3}{4}.
\end{align}
Here $\rho$ is given as in \x{p}.

\begin{theor}[\!]\lb{main}
Assume $({\rm H}_1)$ and $({\rm H}_2)$ hold. Furthermore{\rm ,} it
is assumed that for all $i=1,2,\dots,n,$ the following hypotheses
hold{\rm :}

\begin{enumerate}
\leftskip .7pc
\item[$({\rm D}_1)$] There exist a constant
$\alpha>0$ and a continuous function $\psi_{i}\hbox{\rm :}\
\big[\frac{1}{4},\frac{3}{4}\big]\rightarrow (0,\oo)$ such that
\begin{equation*}
\hskip -1.25pc f^{i}(t,{\rm{\bf{u}}})\geq (\rho\alpha)^{p-1}
\psi_{i}(t),\qquad \frac{1}{4}\leq t\leq \frac{3}{4}
\end{equation*}
for all $0\leq u_{j}\leq\alpha~(j\in\{1,2,\dots,n\}\setminus
\{i\})$ and $\rho\alpha\leq u_{i}\leq\alpha,$ and
\begin{equation*}
\hskip -1.25pc \inf_{\frac{1}{4}\leq t\leq
\frac{3}{4}}\gamma_{\psi_{i}}(t)\geq 1.
\end{equation*}

\item[$({\rm D}_2)$] There exist a constant
$\beta>0,~\beta\neq\alpha$ and a continuous function
$\varphi_{i}\hbox{\rm :}\ [0,1]\rightarrow (0,\oo)$ such that
\begin{equation*}
\hskip -1.25pc f^{i}(t,{\rm{\textbf{u}}})\leq
\beta^{p-1}\varphi_{i}(t)\quad \hbox{for}\,\,0\leq t\leq 1
\,\,\hbox{and}\,\,~0<u_{j}\leq \beta,~~j \in \{1,\dots,n\}
\end{equation*}
and
\begin{equation*}
\hskip -1.25pc \int_{0}^{1}\frac{1}{p(s)} \phi^{-1}
\left(\frac{1}{q(s)} \int_{0}^{1}\varphi_{i}(\tau){\rm
d}\tau\right){\rm d}s\leq 1.
\end{equation*}
Then problems {\rm \x{a}} and {\rm \x{b}} have at least one
positive solution $\rm{\bf{u}}$ satisfying
\begin{equation*}
\hskip -1.25pc
\min\{\alpha,\beta\}\leq\|\rm{\bf{u}}\|\leq\max\{\alpha,\beta\}.
\end{equation*}
\end{enumerate}
\end{theor}

\begin{proof}
We assume that $\alpha<\beta.$ The case $\alpha>\beta$ is
analogous.

Define the sets
\begin{equation*}
\Omega^{1}=\{x\in X\hbox{\rm :}\ \|x\|<\alpha\}\quad
\mbox{and}\quad \Omega^{2}=\{x\in X\hbox{\rm :}\ \|x\|<\beta\}.
\end{equation*}

We claim that
\begin{equation*}
\hskip -4pc {\rm{(i)}}\hskip 3pc\|T \textbf{u}\| \geq
\|\textbf{u}\|~\mbox{for}~ \textbf{u}\in \partial_{K}\Omega^{1};
\end{equation*}
and
\begin{equation*}
\hskip -4pc {\rm{(ii)}}\hskip 2.8pc\|T \textbf{u}\| \leq
\|\textbf{u}\|~ \mbox{for}~\textbf{u}\in \partial_{K}\Omega^{2}.
\end{equation*}
First we shall prove (i). Note, from the definition of
$T\textbf{u}$, that $T^{i}\textbf{u}(\sigma_{i})$ is the maximum
of $T^{i}\textbf{u}$ on $[0,1]$. Since $\|\textbf{u}\|=\alpha$,
there exists $i\in\{1,\dots,n\}$ such that $\rho\alpha\leq
u_{i}(t)\leq\alpha,~\frac{1}{4}\leq t\leq \frac{3}{4}$ and $0\leq
u_{j}\leq\alpha~(j\in\{1,2,\dots,n\}\setminus \{i\})$. We consider
three cases.

\begin{case}{\rm
$\sigma_{i}\in\big[\frac{1}{4},\frac{3}{4}\big].$ Then
\begin{align*}
\sup_{0\leq t\leq 1}(T^{i}\textbf{u})(t) &\ge
\frac{1}{2}\left[\int_{\frac{1}{4}}^{\sigma_{i}}\frac{1}{p(s)}\phi^{-1}
\left(\frac{1}{q(s)}\int_{s}^{\sigma_{i}}f^{i}(\tau,
{\bf{u}}(\tau)){\rm d}\tau\right){\rm d}s \right.\\[.5pc]
&\quad\, \left. + \int^{\frac{3}{4}}_{\sigma_{i}}\frac{1}{p(s)}
\phi^{-1} \left(\frac{1}{q(s)}\int^{s}_{\sigma_{i}}f^{i}
(\tau,{\bf{u}}(\tau)){\rm d}\tau\right){\rm
d}s\right]\nn\\[.5pc]
&\ge \frac{\rho\alpha}{2} \left[\int_{\frac{1}{4}}^{\sigma_{i}}
\frac{1}{p(s)}\phi^{-1} \left(\frac{1}{q(s)}\int_{s}^{\sigma_{i}}
\psi_{i}(\tau){\rm d}\tau\right){\rm
d}s \right.\\[.5pc]
&\quad\, \left. +\int^{\frac{3}{4}}_{\sigma_{i}}\frac{1}{p(s)}
\phi^{-1}\left(\frac{1}{q(s)}\int^{s}_{\sigma_{i}}\psi_{i}(\tau){\rm
d}\tau\right){\rm
d}s\right]\nn\\[.5pc]
&= \alpha\gamma_{\psi_{i}}(\sigma_{i})\geq \alpha
\min_{\frac{1}{4} \leq t\leq
\frac{3}{4}}\gamma_{\psi_{i}}(t)\geq\alpha.\nn
\end{align*}}
\end{case}

\begin{case}{\rm
$\sigma_{i}>\frac{3}{4}.$ Then
\begin{align}
\sup_{0\leq t\leq 1}(T^{i}{\bf{u}})(t) &\ge
\int_{\frac{1}{4}}^{\frac{3}{4}}\frac{1}{p(s)}\phi^{-1} \left(
\frac{1}{q(s)}\int_{s}^{\frac{3}{4}}f^{i}(\tau,
{\bf{u}}(\tau)){\rm d}\tau\right){\rm d}s\nn\\[.5pc]
&\ge \rho\alpha\int_{\frac{1}{4}}^{\frac{3}{4}}
\frac{1}{p(s)}\phi^{-1} \left(\frac{1}{q(s)}
\int_{s}^{\frac{3}{4}}\psi_{i}(\tau){\rm d}\tau \right){\rm d}s\nn\\[.5pc]
&= 2\alpha\gamma_{\psi_{i}}\left(\frac{3}{4}
\right)\geq\alpha\min_{\frac{1}{4}\leq t\leq
\frac{3}{4}}\gamma_{\psi_{i}}(t)\geq\alpha.\nn
\end{align}}
\end{case}

\begin{case}{\rm
$\sigma_{i}<\frac{1}{4}.$ Then
\begin{align}
\sup_{0\leq t\leq 1}(T^{i} {\bf{u}})(t) &\ge
\int_{\frac{1}{4}}^{\frac{3}{4}} \frac{1}{p(s)}\phi^{-1}
\left(\frac{1}{q(s)} \int^{s}_{\frac{1}{4}}f^{i}
(\tau,{\bf{u}}(\tau)){\rm d}\tau\right){\rm d}s\nn\\[.5pc]
&\ge \rho\alpha\int_{\frac{1}{4}}^{\frac{3}{4}}
\frac{1}{p(s)}\phi^{-1} \left(\frac{1}{q(s)}\int^{s}_{\frac{1}{4}}
\psi_{i}(\tau){\rm d}\tau\right){\rm d}s\nn\\[.5pc]
&= 2\alpha\gamma_{\psi_{i}}
\left(\frac{1}{4}\right)\geq\alpha\min_{\frac{1}{4}\leq t\leq
\frac{3}{4}}\gamma_{\psi_{i}}(t)\geq\alpha.\nn
\end{align}}
\end{case}

Hence $\sup_{0\leq t\leq 1}(T^{i}\textbf{u})(t)\geq\alpha.$ This
implies that (i) holds.

Next we prove (ii). In fact, for any $\textbf{u}\in
\partial_{K}\Omega^{2},$ we have $|u_{i}|_{0}\leq\beta$ for each $i\in
\{1,\dots,n\}$. Fix $\,i\in \{1,\dots,n\}$. Then
\begin{align}
(T^{i}\textbf{u}) (t) &\le \int_{0}^{1}\frac{1}{p(s)}\phi^{-1}
\left(\frac{1}{q(s)}\int_{0}^{1}f^{i}(\tau,\textbf{u}(\tau))
{\rm d}\tau\right){\rm d}s\nn
\end{align}
\begin{align}
&\le \beta\int_{0}^{1}\frac{1}{p(s)}\phi^{-1} \left(
\frac{1}{q(s)}
\int_{0}^{1}\varphi_{i}(\tau) {\rm d}\tau\right){\rm d}s\nn\\[.5pc]
&\le \beta.\nn
\end{align}
Therefore, $|T_{i}\textbf{u}|_{0}\leq\|\textbf{u}\|$ for each
$i=1,2,\dots, n.$ This implies that (ii) holds.

Now Theorem~\ref{cone} guarantees that $T$ has a fixed point
$\textbf{u}\in\overline{\Omega^2}_{K} \setminus \Omega^1_{K}.$
Thus $\alpha\leq\|\textbf{u}\|\leq \beta$. Clearly, $\textbf{u}$
is a positive solution of {\rm \x{a} and \x{b}}.\hfill $\Box$
\end{proof}

Next we employ Theorem~\ref{main} to establish the existence of
positive solutions of the following problem:
\begin{equation}\lb{n}
\left\{\begin{array}{ll}
(q(t)\phi(p(t)u'_{1}(t)))'+h_{1}(t)g^{1}(\textbf{u})=0,\quad
0<t<1,\\[.1pc]
\dots\\[.2pc]
(q(t)\phi(p(t)u'_{n}(t)))'+h_{n}(t)g^{n}(\textbf{u})=0,\quad
0<t<1.\end{array}\right.
\end{equation}
We assume the following conditions:

\begin{itemize}
\leftskip .2pc
\item[$({\textmd{H}_{3}})$] $g^{i}\hbox{\rm :}\ \R^{n}_{+}\rightarrow\R_{+}$ is
continuous with $g^{i}(\textbf{u})>0$ for $\|\textbf{u}\|>0,
i=1,2,\dots,n.$

\item[$({\textmd{H}_{4}})$] $h_{i}(t)\hbox{\rm :}\ [0,1]\rightarrow\R_{+}$ is
continuous and $h_{i}(t)\not\equiv 0$ on any subinterval of [0,1],
$i=1,2,\dots,n.$
\end{itemize}
In order to state the result, we introduce the notation
\begin{align*}
g_{0}^{i} &= \lim_{\|\textbf{u}\|\rightarrow 0}
\frac{g^{i}(\textbf{u})}{\phi(\|\textbf{u}\|)},\quad
g_{\oo}^{i}=\lim_{\|\textbf{u}\|\rightarrow
\oo}\frac{g^{i}(\textbf{u})}{\phi(\|\textbf{u}\|)},~\textbf{u}\in
\R^{n}_{+},~i=1, 2, \dots, n.\\[.5pc]
A_{i} &= \int_{0}^{1}\frac{1}{p(s)}
\phi^{-1}\left(\frac{1}{q(s)}\int_{0}^{1}h_{i}(\tau){\rm
d}\tau\right){\rm
d}s,\quad B_{i}=\min_{\frac{1}{4}\leq t\leq \frac{3}{4}}\gamma_{h_{i}}(t)\\[.5pc]
A &= \max\{A_{i}, i=1,2,\dots, n\},\quad B=\min\{B_{i},
i=1,2,\dots, n\}.
\end{align*}

\begin{theor}[\!]\lb{sep} Suppose that conditions
$({\rm H}_{2})$ and $({\rm H}_{4})$ hold. Then problems~{\rm
\x{n}} and {\rm \x{b}} have at least one positive solution
$\rm{\textbf{u}}$ with ${\rm{\textbf{u}}}(t)\not\equiv 0$ for
$t\in(0,1)$ if one of the following conditions holds.

\begin{itemize}
\leftskip 1.8pc
\item[$({\rm h}_{1})$] $0\leq
g^{i}_{0}<\big(\frac{1}{A_{i}}\big)^{p-1}$ and
$\big(\frac{1}{B_{i}}\big)^{p-1}<g^{i}_{\oo}\leq \oo,~i=1,2,\dots,
n${\rm ;}

\item[$({\rm h}_{2})$] $0\leq
g^{i}_{\oo}<\big(\frac{1}{A_{i}}\big)^{p-1}$ and
$\big(\frac{1}{B_{i}}\big)^{p-1}<g^{i}_{0}\leq \oo,~i=1,2,\dots,
n$.
\end{itemize}
\end{theor}

\begin{proof} To see this, we will apply Theorem~\ref{main} with
$f^{i}(t,\textbf{u})=h_{i}(t)g^{i}(\textbf{u}),~i=1,2,\dots, n.$
We assume that $({\rm h}_{1})$ holds. The case when $({\rm
h}_{2})$ holds is similar.

From the first part of $({\rm h}_{1})$, there exists $\beta>0$
such that $g^{i}(\textbf{u})\leq
\big(\frac{1}{A_{i}}\big)^{p-1}\beta$ for
$\|\textbf{u}\|\leq\beta.$ Choose
$\varphi_{i}(t)=\big(\frac{1}{A_{i}}\big)^{p-1}h_{i}(t)$ for
$i=1,2,\dots, n.$ Fix $\,i\in \{1,\dots, n\}$. Then
\begin{align*}
f^{i}(t,\textbf{u}) &= h_{i}(t)g^{i}(\textbf{u})\leq
\left(\frac{1}{A_{i}}\right)^{p-1}\beta h_{i}(t)=\beta
\varphi_{i}(t)\\[.5pc]
&\quad\, \hbox{ if }\,\, 0\leq t\leq 1 \,\,\hbox{ and
}\,\,~0<u_{j}\leq\beta \,\,\hbox{ for }\,\,j\in \{1,\dots,n\}
\end{align*}
and
\begin{align*}
&\int_{0}^{1}\frac{1}{p(s)} \phi^{-1}\left(\frac{1}{q(s)}
\int_{0}^{1}\varphi_{i}(\tau){\rm d}\tau\right){\rm
d}s\\[.5pc]
&\quad\, =\frac{1}{A_{i}} \int_{0}^{1}\frac{1}{p(s)}
\phi^{-1}\left(\frac{1}{q(s)} \int_{0}^{1}h_{i}(\tau){\rm
d}\tau\right){\rm d}s=1.
\end{align*}
Thus hypothesis $({\rm D}_{2})$ holds.

From the second part of $({\rm h}_{1})$, there exists $\alpha>0$
such that $\alpha>\max\{\rho^{-1}\beta,\beta\}$ and
$g^{i}(\textbf{u})\geq\big(\frac{1}{B_{i}}\big)^{p-1}(\rho\alpha)^{p-1}$
for $\|\textbf{u}\|\geq\rho\alpha,~i=1, 2, \dots, n.$

Thus $g^{i}(\textbf{u})\geq \big(\frac{1}{B_{i}}\big)^{p-1}
(\rho\alpha)^{p-1}$ for $u_{i}\geq\rho\alpha,~i=1, 2, \dots, n.$
Choose $\psi_{i}(t)=\big(\frac{1}{{B_{i}}}\big)^{p-1}h_{i}(t)$ for
$i=1,2,\dots, n.$ Then
\begin{align*}
f^{i}(t,\textbf{u}) &= h_{i}(t)g^{i}(\textbf{u})\geq
\left(\frac{1}{{B_{i}}} \right)^{p-1}(\rho\alpha)^{p-1}
h_{i}(t)=(\rho\alpha)^{p-1} \psi_{i}(t),\\[.5pc]
&\quad\, \frac{1}{4}\leq t\leq \frac{3}{4},~u_{i}\geq\rho\alpha,
\end{align*}
(so~in~particular~for~$\rho\alpha\leq u_{i}\leq\alpha)$ and for
$\frac{1}{4}\leq t\leq \frac{3}{4}$,
\begin{align}
\gamma_{\psi_{i}}(t) &= \frac{\rho}{2}
\left[\int_{\frac{1}{4}}^{t} \frac{1}{p(s)}\phi^{-1}
\left(\frac{1}{q(s)} \int_{s}^{t}\psi_{i}(\tau){\rm
d}\tau\right){\rm d}s\right.\nn\\[.5pc]
&\quad\, \left. + \int^{\frac{3}{4}}_{t}\frac{1}{p(s)}
\phi^{-1}\left(\frac{1}{q(s)}\int^{s}_{t}\psi_{i}(\tau){\rm d}
\tau\right){\rm d}s\right]\nn\\[.5pc]
&= \frac{1}{B_{i}}\gamma_{h_{i}}(t).\nn
\end{align}
Thus
\begin{equation*}
\inf_{\frac{1}{4}\leq t\leq \frac{3}{4}} \gamma_{\psi_{i}}(t) =
\frac{1}{B_{i}}\inf_{\frac{1}{4}\leq t\leq \frac{3}{4}}
\gamma_{h_{i}}(t)=1.
\end{equation*}
This implies that hypothesis $({\rm D}_{1})$ holds. The result now
follows from Theorem~\ref{main}.\hfill $\Box$
\end{proof}

Next we consider the nonlinear eigenvalue problems \x{i} and
\x{b}. By applying Theorem~\ref{sep}, we easily get the following
result.

\begin{theor}[\!]\lb{eig} Suppose that conditions
$({\rm H}_{2})$ and $({\rm H}_{4})$ hold. Then problems~{\rm
\x{i}} and {\rm \x{b}} have at least one positive solution for
each
\begin{equation}
\lb{s}\lambda\in\left(\disp\frac{1}{B^{p-1}\disp
\min_{i=1,2,\dots,n}\{g^{i}_{\oo}\}},\disp\frac{1}
{A^{p-1}\disp\max_{i=1,2,\dots,n}\{g_{0}^{i}\}}\right)
\end{equation}
if ${1}/({B^{p-1}\min_{i=1,2,\dots,n}\{g^{i}_{\oo}\}})<
{1}/({A^{p-1}\max_{i=1,2,\dots,n}\{g_{0}^{i}\}})$. The same result
remains valid for each
\begin{equation}\lb{t}
\lambda\in\left(\disp\frac{1}{B^{p-1}\disp\min_{i=1,2,\dots,n}
\{g^{i}_{0}\}},\disp\frac{1}
{A^{p-1}\disp\max_{i=1,2,\dots,n}\{g_{\oo}^{i}\}}\right)
\end{equation}
if ${1}/({B^{p-1} \min_{i=1,2,\dots,n} \{g^{i}_{0}\}})< {1}/
({A^{p-1}\max_{i=1,2,\dots,n} \{g_{\oo}^{i}\}})$. Here we write
$1/g^{i}_{\alpha}=0$ if $g^{i}_{\alpha}=\oo$ and
$1/g^{i}_{\alpha}=\oo$ if $g^{i}_{\alpha}=0,$ where
$\alpha=0,\oo.$
\end{theor}

\begin{proof}We consider the case \x{s}. The case \x{t} is similar. If
$\lambda$ satisfies \x{s}, then
\begin{equation*}
\lambda g_{0}^{i}\leq\lambda\max_{i=1,2,\dots,n}
\{g_{0}^{i}\}<\frac{1}{A^{p-1}}
\leq\left(\frac{1}{A_{i}}\right)^{p-1},\qquad i=1,2,\dots,n,
\end{equation*}
and
\begin{equation*}
\lambda g_{\oo}^{i}\geq
\lambda\min_{i=1,2,\dots,n}\{g_{\oo}^{i}\}>\frac{1}{B^{p-1}} \geq
\left(\frac{1}{B_{i}}\right)^{p-1},\qquad i=1,2,\dots,n.
\end{equation*}
So Theorem \ref{sep} applies directly.\hfill $\Box$
\end{proof}

\begin{coro}$\left.\right.$\vspace{.5pc}

\noindent Problems {\rm\x{i}} and {\rm \x{b}} have at least one
positive solution for each $\lambda\in(0,\oo)$ if one of the
following two conditions holds{\rm :}

\begin{enumerate}
\renewcommand\labelenumi{\rm (\roman{enumi})}
\leftskip .15pc
\item $g_{\oo}^{i}=\oo$ and $g_{0}^{i}=0,$
$i=1,2,\dots,n;$

\item $g_{0}^{i}=\oo$ and $g_{\oo}^{i}=0,$
$i=1,2,\dots,n.$\vspace{-.5pc}
\end{enumerate}
\end{coro}

Finally, it is worth remarking here that, we can apply
Theorems~\ref{sep} and \ref{eig} to study the existence of
positive radial solutions for the nonlinear elliptic systems \x{c}
and \x{d}.

In fact, a radial solution of \x{c} and \x{d} can be considered as
a solution of the system
\begin{equation}\lb{e}
\left\{\begin{array}{ll} (r^{N-1}\phi(u'_{1}
(r)))'+r^{N-1}k_{1}(r)g^{1}(\textbf{u})=0,\\[.1pc]
\dots\\[.2pc]
(r^{N-1}\phi(u'_{n}
(r)))'+r^{N-1}k_{n}(r)g^{n}(\textbf{u})=0,\end{array}\right.
\end{equation}
$0<R_{1}<r<R_{2}<\oo,$ with the following boundary condition
\begin{equation}\lb{f}
\textbf{u}(R_{1})=0,\quad \textbf{u}(R_{2})=0.
\end{equation}
Applying the change of variables, $r=(R_{2}-R_{1})t+R_{1},$ we can
transform \x{e} and \x{f} into the form
\begin{equation}\lb{g}
\left\{\begin{array}{ll} (q(t)\phi(\xi u'_{1}(t)))'+
h_{1}(t)g^{1}(\textbf{u})=0,\quad 0<t<1,\\[.1pc]
\dots\\[.2pc]
(q(t)\phi(\xi u'_{n}(t)))'+h_{n}(t)g^{n}(\textbf{u})=0,\quad
0<t<1,\end{array}\right.
\end{equation}
with boundary condition \x{b}, where
\begin{equation*}
q(t)=((R_{2}-R_{1})t+R_{1})^{N-1},~\zeta=\frac{1}{R_{2}-R_{1}}
\end{equation*}
and
\begin{equation*}
\hskip -4pc
h_{i}(t)=(R_{2}-R_{1})((R_{2}-R_{1})t+R_{1})^{N-1}k_{i}((R_{2}-
R_{1})t+R_{1}),~i=1,2,\dots,n,
\end{equation*}
which correspond to problems \x{a} and \x{b} with $p(t)=\xi$ and
$f^{i}(t,\textbf{u})=h_{i}(t)g^{i}(\textbf{u}),~i=1,2,\dots,n.$

Thus, we only need to consider problems~\x{g} and \x{b} and we can
obtain results similar to those in Theorems~\ref{sep} and
\ref{eig}. Here we omit the details.

\section*{Acknowledgement}

This work is supported by the National Natural Science Foundation
of China (Item no.~10325102 and no.~10531010).

\end{document}